\documentclass[]{amsart}
\usepackage{amsmath}
\usepackage{amsfonts,amssymb,latexsym}
\usepackage{amscd}
\usepackage{mathrsfs}
\usepackage{amsthm}
\theoremstyle{plain}
\newtheorem{theorem}{\rmfamily\bfseries{Theorem}}[section] 
\newtheorem{corollary}[theorem]{\rmfamily\bfseries{Corollary}} 

\newtheorem{lemma}[theorem]{\rmfamily\bfseries{Lemma}} 
\newtheorem{proposition}[theorem]{\rmfamily\bfseries{Proposition}}

\newtheorem{definition}[theorem]{\rmfamily\bfseries{Definition}}%[section]

\theoremstyle{remark}

\newtheorem{remark}{Remark}

\usepackage{graphics}
\numberwithin{equation}{section}

\newcommand{\mona}{\mbox{\LARGE\itshape a}}

\def\un{1\kern-3pt \rm I}
\def\ptoday{{\ifcase\month 
\or January, \or February, \or March, \or April,\or May, 
\or June, \or July, \or August, \or September, \or October, 
\or November, \or December,\fi\ \number \year}}
\newcommand{\oN}{{\mathbb N}}
\newcommand{\oR}{{\mathbb R}}

\newcommand{\oC}{{\mathbb C}}

\newcommand{\supp}{\mathrm{supp}}
   
\begin{document}
%%%%%%%%%%%%%%%%%%%%%%%%%%%%%%%%%%%%%%%%%%%%%%%%%%%%%%%%%%%%%%%%%%%%%%%%%%%%%%%%
%{\hfill
%\parbox{50mm}{{\sf CEFT-SFM-DHTF06/1}\\{\sf Revised version}} \vspace{8mm}}

\title[Holomorphic Extension Theorem for Ultrahyperfunctions]
      {{Holomorphic Extension Theorem for\\ Tempered Ultrahyperfunctions}}
        
\author{Daniel H.T. Franco}
\address{Universidade Federal de Vi\c cosa\\
         Departamento de F\'\i sica\\
         Avenida Peter Henry Rolfs s/n\\  
         Campus Universit\'ario, Vi\c cosa, MG, Brasil, CEP:36570-000.}
\email{dhtfranco@gmail.com}

\date{\today}
\thanks{This work was supported by the Funda\c c\~ao de Amparo \`a Pesquisa do
Estado de Minas Gerais (FAPEMIG) agency, grant CEX00012/07.}

\begin{abstract}
In this paper we are concerned with the space of tempered ultrahyperfunctions
corresponding to a proper open convex cone. A holomorphic extension theorem (the version
of the celebrated edge of the wedge theorem) will be given for this setting. As
application, a version is also given of the principle of determination of an analytic
function by its values on a non-empty open real set. The paper finishes with the
generalization of holomorphic extension theorem \`a la Martineau.
\end{abstract}

\maketitle

%\centerline{\underline{\bf DRAFT VERSION}}

%%%%%%%%%%%%%%%%%%%%%%%%%%%%%%%%%%%%%%%%%%%%%%%%%%%%%%%
\section{Introduction}
%\hspace\parindent
%%%%%%%%%%%%%%%%%%%%%%%%%%%%%%%%%%%%%%%%%%%%%%%%%%%%%%%
Sebasti\~ao e Silva~\cite{Tiao1}, \cite{Tiao2} and Hasumi~\cite{Hasumi} introduced the
space of tempered ultradistributions, which has been studied by many authors,
among others we refer the reader to~\cite{Zie}-\cite{DanHenri}. Here, as
Morimoto~\cite{Mari1} and \cite{Mari3}, we shall refer to the tempered ultradistributions
as {\it tempered ultrahyperfunctions} in order to distinguish them from various
other classes of ultradistributions which have been described as tempered (see,
for example, Pathak~\cite{Pathak} and Pilipovic~\cite{Pili}). Tempered ultrahyperfunctions
are the strong dual of the space of test functions of rapidly decreasing entire
functions in any horizontal strip. While Sebasti\~ao e Silva~\cite{Tiao1} used extension
procedures for the Fourier transform combined with holomorphic representations and
considered the 1-dimensional case, Hasumi~\cite{Hasumi} used duality arguments in
order to extend the notion of tempered ultrahyperfunctions for the case of $n$ dimensions
(see also~\cite[Section 11]{Tiao2}). In a brief tour, Marimoto~\cite{Mari1}
gave some more precise informations concerning the work of Hasumi.
More recently, the relation between the tempered ultrahyperfunctions and
Schwartz distributions and some major results, as the kernel theorem and the
Fourier-Laplace transform have been established by Br\"uning-Nagamachi
in~\cite{BruNa1}. Earlier, some precisions on the Fourier-Laplace transform
theorem for tempered ultrahyperfunctions were given by Carmichael~\cite{Carmi1}
(see also~\cite{DanHenri}), by considering the theorem in its simplest form,
{\em i.e.}, the equi\-va\-lence between support properties of a distribution
in a closed convex cone and the holomorphy of its Fourier-Laplace transform
in a suitable tube with conical basis. In this more general setting, which
includes the results of Sebasti\~ao e Silva and Hasumi as special cases, Carmichael
obtained new representations of tempered ultrahyperfunctions which were not considered
in Refs.~\cite{Tiao1}, \cite{Tiao2}, \cite{Hasumi}. 

The purpose of this paper is to prove the edge of the wedge theorem for the setting
of tempered ultrahyperfunctions corresponding to a proper open convex cone. This classical
theorem in complex analysis, discovered by theoretical physicists in 1950's~\cite{EWT},
deals with the question about the principle of holomorphic continuation of functions of
several complex variables, which arose in physics in the study of the Wightman functions
and Green functions, or in connection with the dispersion relations in quantum field
theory. It should be mentioned that other versions of the theorem for tempered
ultrahyperfunctions can be found in Refs.~\cite{Mari3}, \cite{Suwa}. Our approach to
this problem is different from that taken in Refs.~\cite{Mari3}, \cite{Suwa}. Our
construction parallels that of Carmichael~\cite{Carmi1}, \cite{Carmi2}, \cite{Carmi3}
and, in particular, the proof of the edge of the wedge theorem is inspired by Carmichael's work~\cite{Carmi3}. As an immediate application of the edge of the wedge theorem, we give
also a proof of the principle of determination of an analytic function by its values on
a non-empty open real set. We finish with a generalized version of holomorphic extension
theorem \`a la Martineau.

We note that the results obtained here are of interest in the construction
and study of {\em quasilocal} quantum field theories (where the fields are
localizable only in regions greater than a certain scale of nonlocality), since
Br\"uning-Nagamachi~\cite{BruNa1} have recently shown the importance of
tempered ultrahyperfunctions for quantum field theories with a fundamental length.
This is the case of a quantum field theory in non-commutative spacetimes~\cite{DZH}.

%%%%%%%%%%%%%%%%%%%%%%%%%%%%%%%%%%%%%%%%%%%%%%%%%%%%%%%%%%%%%%%%%%%%%%%%%%%%%
\section{Notation and Definitions}
\label{SecUltra}
%\hspace*{\parindent}
%%%%%%%%%%%%%%%%%%%%%%%%%%%%%%%%%%%%%%%%%%%%%%%%%%%%%%%%%%%%%%%%%%%%%%%%%%%%%
The following multi-index notation is used without further explanation. Let
$\oR^n$ (resp. $\oC^n=\oR^n+i\oR^n$) be the real (resp. complex) $n$-space whose generic
points are denoted by $x=(x_1,\ldots,x_n)$ (resp. $z=(z_1,\ldots,z_n)$), such that
$x+y=(x_1+y_1,\ldots,x_n+y_n)$, $\lambda x=(\lambda x_1,\ldots,\lambda x_n)$,
$x \geq 0$ means $x_1 \geq 0,\ldots,x_n \geq 0$, $\langle x,y \rangle=x_1y_1+\cdots+x_ny_n$
and $|x|=|x_1|+\cdots+|x_n|$. Moreover, we define
$\alpha=(\alpha_1,\ldots,\alpha_n) \in \oN^n_o$, where $\oN_o$ is the set
of non-negative integers, such that the length of $\alpha$ is the corresponding
$\ell^1$-norm $|\alpha|=\alpha_1+\cdots +\alpha_n$, $\alpha+\beta$ denotes
$(\alpha_1+\beta_1,\ldots,\alpha_n+\beta_n)$, $\alpha \geq \beta$ means
$(\alpha_1 \geq \beta_1,\ldots,\alpha_n \geq \beta_n)$, $\alpha!=
\alpha_1! \cdots \alpha_n!$, $x^\alpha=x_1^{\alpha_1}\ldots x_n^{\alpha_n}$,
and
\[
D^\alpha \varphi(x)=\frac{\partial^{|\alpha|}\varphi(x_1,\ldots,x_n)}
{\partial x_1^{\alpha_1}\partial x_2^{\alpha_1}\ldots\partial x_n^{\alpha_n}}\,\,.
\]
Let $\Omega$ be a set in $\oR^n$. Then we denote by $\Omega^\circ$ the interior
of $\Omega$ and by $\overline{\Omega}$ the closure of $\Omega$. For $r > 0$, we
denote by $B(x_o;r)=\bigl\{x \in \oR^n \mid |x-x_o| < r\bigr\}$ an open ball
and by $B[x_o;r]=\bigl\{x \in \oR^n \mid |x-x_o| \leq r\bigr\}$ a closed ball,
with center at point $x_o$ and of radius $r$, respectively.

We consider two $n$-dimensional spaces -- $x$-space and $\xi$-space -- with the
Fourier transform defined
\[
\widehat{f}(\xi)={\mathscr F}[f(x)](\xi)=
\int_{\oR^n} f(x)e^{i \langle \xi,x \rangle} d^nx\,\,,
\]
while the Fourier inversion formula is
\[
f(x)={\mathscr F}^{-1}[\widehat{f}(\xi)](x)= \frac{1}{(2\pi)^n}
\int_{\oR^n} \widehat{f}(\xi)e^{-i \langle \xi,x \rangle} d^n\xi\,\,.
\]
The variable $\xi$ will always be taken real while $x$ will also be
complexified -- when it is complex, it will be noted $z=x+iy$. The
above formulas, in which we employ the symbolic ``function notation,''
are to be understood in the sense of distribution theory.

%%%%%%%%%%%%%%%%%%%%%%%%%%%%%%%%%%%%%%%%%%%%%%%%%%%%%%%%%%%%%%%%%%%%%%%%%%%%%
\section{Tempered Ultrahyperfunctions}
\label{Sec3}
%\hspace*{\parindent}
%%%%%%%%%%%%%%%%%%%%%%%%%%%%%%%%%%%%%%%%%%%%%%%%%%%%%%%%%%%%%%%%%%%%%%%%%%%%%
Since the theory of ultrahyperfunctions is not too well known, we shall introduce
briefly in this section some definitions and basic properties of the
tempered ultrahyperfunction space of Sebasti\~ao e Silva~\cite{Tiao1}, \cite{Tiao2}
and Hasumi~\cite{Hasumi} (we indicate the Refs. for more details) used
throughout the paper.
To begin with, we shall consider the function
\[
h_{K}(\xi)=\sup_{x \in K} |\langle \xi,x \rangle|\,\,,
\quad \xi \in \oR^n\,\,,
\]
where $K$ is a compact set in $\oR^n$. One calls $h_{K}(\xi)$ the {\it supporting
function} of $K$. We note that $h_{K}(\xi) < \infty$ for every $\xi \in \oR^n$ since
$K$ is bounded. For sets $K=\bigl[-k,k\bigr]^n$, $0 < k < \infty$, the supporting
function $h_{K}(\xi)$ can be easily determined:
\[
h_{K}(\xi)=\sup_{x \in K} |\langle \xi,x \rangle|=
k|\xi|\,\,,\quad \xi \in \oR^n\,\,,\quad |\xi|=\sum_{i=1}^n|\xi_i|\,\,.
\]
Let $K$ be a convex compact subset of $\oR^n$,
then $H_b(\oR^n;K)$ ($b$ stands for bounded) defines the space of all
functions $\in C^\infty(\oR^n)$ such that $e^{h_K(\xi)}D^\alpha \varphi(\xi)$
is bounded in $\oR^n$ for any multi-index $\alpha$. One defines in
$H_b(\oR^n;K)$ seminorms
\begin{equation}
\|\varphi\|_{K,N}=\sup_{{\substack{\xi \in \oR^n \\ \alpha \leq N}}}
\bigl\{e^{h_K(\xi)}|D^\alpha \varphi(\xi)|\bigr\} < \infty\,\,,
\quad N=0,1,2,\ldots\,\,.
\label{snorma2}
\end{equation}

If $K_1 \subset K_2$ are two compact convex sets, then
$h_{K_1}(\xi) \leq h_{K_2}(\xi)$, and thus the canonical
injection $H_b(\oR^n;K_2) \hookrightarrow H_b(\oR^n;K_1)$
is continuous. Let $O$ be a convex open set of $\oR^n$.
To define the topology of $H(\oR^n;O)$ it suffices to let $K$ range
over an increasing sequence of convex compact subsets $K_1,K_2,\ldots$
contained in $O$ such that for each $i=1,2,\ldots$,
$K_i \subset K_{i+1}^\circ$ and ${O}=\bigcup_{i=1}^\infty K_i$.
Then the space $H(\oR^n;O)$ is the projective limit of the
spaces $H_b(\oR^n;K)$ according to restriction mappings
above, {\em i.e.}
\begin{equation}
H(\oR^n;O)=\underset{K \subset {O}}{\lim {\rm proj}}\,\,
H_b(\oR^n;K)\,\,,
\label{limproj2}
\end{equation} 
where $K$ runs through the convex compact sets contained in $O$.

\begin{theorem}[\cite{Hasumi}, \cite{Mari1}, \cite{BruNa1}]
The space ${\mathscr D}({\oR^n})$ of all $C^\infty$-functions
on $\oR^n$ with compact support is dense in $H(\oR^n;K)$ and $H(\oR^n;O)$.
The space $H(\oR^n;\oR^n)$ is dense in $H(\oR^n;O)$ and in
$H(\oR^n;K)$, and $H(\oR^m;\oR^m) \otimes H(\oR^n;\oR^n)$ is dense in
$H(\oR^{m+n};\oR^{m+n})$.
\label{theoINJ}
\end{theorem}

From Theorem \ref{theoINJ} we have the following injections~\cite{Mari1}:
\[
H^\prime(\oR^n;K) \hookrightarrow H^\prime(\oR^n;\oR^n)
\hookrightarrow {\mathscr D}^\prime(\oR^n)\,\,,
\]
and
\[
H^\prime(\oR^n;O) \hookrightarrow H^\prime(\oR^n;\oR^n)
\hookrightarrow {\mathscr D}^\prime(\oR^n)\,\,.
\]

A distribution $V \in H^\prime(\oR^n;O)$ may be expressed as a finite
order deri\-va\-ti\-ve of a continuous function of exponential growth
\[
V=D^\gamma_\xi[e^{h_K(\xi)}g(\xi)]\,\,,
\]
where $g(\xi)$ is a bounded continuous function. For $V \in
H^\prime(\oR^n;O)$ the follo\-wing result is known:

\begin{lemma}[\cite{Mari1}]
A distribution $V \in {\mathscr D}^\prime(\oR^n)$ belongs to $H^\prime(\oR^n;O)$
if and only if there exists a multi-index $\gamma$, a convex compact set $K \subset O$
and a bounded continuous function $g(\xi)$ such that
\[
V=D^\gamma_\xi[e^{h_K(\xi)}g(\xi)]\,\,.
\]
\label{lemmaMari}
\end{lemma}

In the space $\oC^n$ of $n$ complex variables $z_i=x_i+iy_i$,
$1 \leq i \leq n$, we denote by $T(\Omega)=\oR^n+i\Omega \subset \oC^n$
the tubular set of all points $z$, such that $y_i={\text{Im}}\,z_i$ belongs
to the domain $\Omega$, {\em i.e.}, $\Omega$ is a connected open set in $\oR^n$
called the basis of the tube $T(\Omega)$. Let $K$ be a convex compact
subset of $\oR^n$, then ${\mathfrak H}_b(T(K))$ defines
the space of all continuous functions $\varphi$ on $T(K)$ which are holomorphic
in the interior $T(K^\circ)$ of $T(K)$ such that the estimate
\begin{equation}
|\varphi(z)| \leq {\boldsymbol{\sf M}}_{_{K,N}}(\varphi) (1+|z|)^{-N}
\label{est}
\end{equation}
is valid. The best possible constants in (\ref{est}) are given by a family of seminorms in
${\mathfrak H}_b(T(K))$
\begin{equation}
\|\varphi\|_{K,N}=\inf\Bigl\{{\boldsymbol{\sf M}}_{_{K,N}}(\varphi) \mid \sup_{z \in T(K)}
\bigl\{(1+|z|)^{N}|\varphi(z)|\bigr\} < \infty, N=0,1,2,\ldots \Bigr\}\,\,.
\label{snorma1}
\end{equation}

If $K_1 \subset K_2$ are two convex compact sets, we have that the canonical injection
\begin{equation}
{\mathfrak H}_b(T(K_2)) \hookrightarrow {\mathfrak H}_b(T(K_1))\,\,,
\label{canoinj}
\end{equation}
is continuous.

Given that the spaces ${\mathfrak H}_b(T(K_i))$ are Fr\'echet spaces, with topology
defined by the seminorms (\ref{snorma1}), the space ${\mathfrak H}(T({O}))$ is
characterized as a projective limit of Fr\'echet spaces: 
\begin{equation}
{\mathfrak H}(T({O}))=\underset{K \subset {O}}{\lim {\rm proj}}\,\,
{\mathfrak H}_b(T(K))\,\,,
\label{limproj1}
\end{equation}
where $K$ runs through the convex compact sets contained in $O$ and
the projective limit is taken following the restriction mappings above.

Let $K$ be a convex compact set in $\oR^n$. Then the space ${\mathfrak H}(T(K))$
is characterized as a inductive limit
\begin{equation}
{\mathfrak H}(T(K))=\underset{K_1 \supset K}{\lim {\rm ind}}\,\,
{\mathfrak H}_b(T(K_1))\,\,,
\label{limind1}
\end{equation}
where $K_1$ runs through the convex compact sets such that $K$ is contained
in the interior of $K_1$ and the inductive limit is taken following the restriction
mappings (\ref{canoinj}).

For any element $U \in {\mathfrak H}^\prime$, its Fourier transform is
defined to be a distribution $V$ of exponential growth, such that the
Parseval-type relation
\begin{equation}
\langle V,\varphi \rangle=\langle U,\psi \rangle\,\,,\quad
\varphi \in H\,\,,\,\,\psi={\mathscr F}[\varphi] \in {\mathfrak H}\,\,,
\label{PRel1}
\end{equation}
holds. In the same way, the inverse Fourier transform of a distribution $V$ of
exponential growth is defined by the relation 
\begin{equation}
\langle U,\psi \rangle=\langle V,\varphi \rangle\,\,,\quad
\psi \in {\mathfrak H}\,\,,\,\,\varphi={\mathscr F}^{-1}[\psi] \in H\,\,.
\label{PRel2}
\end{equation}

It follows from the Fourier transform and Theorem \ref{theoINJ} the

\begin{theorem}[\cite{Mari1}, \cite{BruNa1}]
${\mathfrak H}(T(\oR^n))$ is dense in ${\mathfrak H}(T(O))$ and
in ${\mathfrak H}(T(K))$, and ${\mathfrak H}(T(\oR^{m+n}))$ is dense
in ${\mathfrak H}(T(O))$.
\label{theoINJE}
\end{theorem}

\begin{proposition}[\cite{Mari1}]
If $f \in H(\oR^n;O)$, the Fourier transform of $f$ belongs
to the space ${\mathfrak H}(T(O))$, for any open convex
non-empty set $O \subset \oR^n$. By the dual Fourier transform
$H^\prime(\oR^n;O)$ is topologically isomorphic with the space
${\mathfrak H}^\prime(T(-O))$.
\label{Propo1}
\end{proposition}

\begin{definition}
A tempered ultrahyperfunction is a continuous linear functional defined
on the space of test functions ${\mathfrak H}(T(\oR^n))$ of rapidly decreasing
entire functions in any horizontal strip. 
\label{UHF}
\end{definition}

The space of all tempered ultrahyperfunctions is denoted by ${\mathscr U}(\oR^n)$.
As a matter of fact, these objects are equi\-va\-lence classes of holomorphic
functions defined by a certain space of functions which are analytic in the $2^n$
octants in $\oC^n$ and represent a natural generalization of the notion of
hyperfunctions on $\oR^n$, but are {\it non-localizable}. The space
${\mathscr U}(\oR^n)$ is characterized in the following way~\cite{Hasumi}:
Let $\boldsymbol{{\mathscr H}_\omega}$ be the space of all
functions $f(z)$ such that ({\it i}) $f(z)$ is analytic for $\{z \in \oC^n \mid
|{\rm Im}\,z_1| > p, |{\rm Im}\,z_2| > p,\ldots,|{\rm Im}\,z_n| > p\}$,
({\it ii}) $f(z)/z^p$ is bounded continuous  in
$\{z\in \oC^n \mid |{\rm Im}\,z_1| \geqq p,|{\rm Im}\,z_2| \geqq p,
\ldots,|{\rm Im}\,z_n| \geqq p\}$, where $p=0,1,2,\ldots$ depends on $f(z)$
and ({\it iii}) $f(z)$ is bounded by a power of $z$, $|f(z)|\leq
{\boldsymbol{\sf C}}(1+|z|)^N$, where ${\boldsymbol{\sf C}}$ and $N$ depend on
$f(z)$. Define the {\em kernel} of the mapping $f:{\mathfrak H}(T(\oR^n))
\rightarrow \oC$ by $\boldsymbol{\Pi}$, as the set of all $z$-dependent
pseudo-polynomials, $z\in \oC^n$ (a pseudo-polynomial is a
function of $z$ of the form $\sum_s z_j^s G(z_1,...,z_{j-1},z_{j+1},...,z_n)$,
such that $G(z_1,...,z_{j-1},z_{j+1},...,z_n) \in \boldsymbol{{\mathscr H}_\omega}$).
Then, $f(z) \in \boldsymbol{{\mathscr H}_\omega}$ belongs to the kernel
$\boldsymbol{\Pi}$ if and only if $\langle f(z),\psi(x) \rangle=0$,
with $\psi(x) \in {\mathfrak H}(T(\oR^n))$ and $x={\rm Re}\,z$. The
space of tempered ultrahyperfunctions is the quotient space ${\mathscr U}=
\boldsymbol{{\mathscr H}_\omega}/\boldsymbol{\Pi}$. Thus, we have the

\begin{theorem}[Hasumi~\cite{Hasumi}, Proposition 5] The space of
tempered ultrahyperfunctions ${\mathscr U}$ is algebraically isomorphic
to the space of generalized functions ${\mathfrak H}^\prime$.
\label{HasumiTheo}
\end{theorem}

%%%%%%%%%%%%%%%%%%%%%%%%%%%%%%%%%%%%%%%%%%%%%%%%%%%%%%%%%%%%%%%%%%%%%%%%%%%%%
\section{The Space of Holomorphic Functions $\boldsymbol{{\mathscr H}^o_c}$}
\label{SecTheo1}
%\hspace*{\parindent}
%%%%%%%%%%%%%%%%%%%%%%%%%%%%%%%%%%%%%%%%%%%%%%%%%%%%%%%%%%%%%%%%%%%%%%%%%%%%%
We start by introducing some terminology and simple
facts concerning cones. An open set $C \subset \oR^n$ is called a cone if
$\oR_+\!\cdot\!C \subset C$. A cone $C$ is an open connected cone if $C$ is an
open connected set. Moreover, $C$ is called convex if $C+C \subset C$ and {\it proper}
if it contains no any straight line. A cone $C^\prime$ is called compact in $C$ --
we write $C^\prime \Subset C$ -- if the projection ${\sf pr}{\overline C^{\,\prime}}
\overset{\text{def}}{=}{\overline C^{\,\prime}} \cap S^{n-1} \subset
{\sf pr}C\overset{\text{def}}{=}C \cap S^{n-1}$, where $S^{n-1}$ is the unit sphere
in $\oR^n$. Being given a cone $C$ in $y$-space, we associate with $C$ a closed convex
cone $C^*$ in $\xi$-space which is the set $C^*=\bigl\{\xi \in \oR^n \mid \langle \xi,y
\rangle \geq 0,\forall\,\,y \in C \bigr\}$. The cone $C^*$ is called the {\em dual cone}
of $C$. In the sequel, it will be sufficient to assume for our purposes that the open
connected cone $C$ in $\oR^n$ is an open convex cone with vertex at the origin
and proper. By $T(C)$ we will denote the set $\oR^n+iC \subset \oC^n$.
If $C$ is open and connected, $T(C)$ is called the tubular radial domain in
$\oC^n$, while if $C$ is only open $T(C)$ is referred to as a tubular cone.  In the
former case we say that $f(z)$ has a boundary value $U=BV(f(z))$ in ${\mathfrak H}^\prime$
as $y \rightarrow 0$, $y \in C$ or $y \in C^\prime \Subset C$, respectively, if
for all $\psi \in {\mathfrak H}$ the limit
\[
\langle U, \psi \rangle=\lim_{{\substack{y \rightarrow 0 \\
y \in C~{\rm or}~C^\prime}}} \int_{\oR^n} f(x+iy)\psi(x) d^nx\,\,,
\]
exists. An important example of tubular radial domain used in quantum field theory
is the tubular radial domain with the forward light-cone, $V_+$, as its basis
\[
V_+=\Bigl\{z \in \oC^n \mid {\rm Im}\,z_1 >
\Bigl(\sum_{i=2}^n {\rm Im}^2\,z_i \Bigr)^{\frac{1}{2}}, {\rm Im}\,z_1 > 0 \Bigr\}\,\,.  
\]

We will deal with tubes defined as the set of all points $z \in \oC^n$
such that
\[
T(C)=\Bigl\{x+iy \in \oC^n \mid x \in \oR^n, y \in C, |y| < \delta \Bigr\}\,\,,
\]
where $\delta > 0$ is an arbitrary number. 

Let $C$ be a proper open convex cone, and let $C^\prime \Subset C$.
Let $B[0;r]$ denote a {\bf closed} ball of the
origin in $\oR^n$ of radius $r$, where $r$ is an arbitrary positive
real number. Denote $T(C^\prime;r)=\oR^n+i\bigl(C^\prime \setminus
\bigl(C^\prime \cap B[0;r]\bigr)\bigr)$.
We are going to introduce a space of holomorphic functions
which satisfy certain estimate according to Carmichael~\cite{Carmi1}.
We want to consider the space consisting of holomorphic functions $f(z)$
such that
\begin{equation}
\bigl|f(z)\bigr|\leq {\boldsymbol{\sf M}}(C^\prime)(1+|z|)^N e^{h_{C^*}(y)}
\,\,,\quad z \in T(C^\prime;r)\,\,,
\label{eq31} 
\end{equation}
where $h_{C^*}(y)=\sup_{\xi \in C^*}\langle \xi,y \rangle$ is the supporting
function of $C^*$, ${\boldsymbol{\sf M}}(C^\prime)$ is a constant that depends
on an arbitrary compact cone $C^\prime$ and $N$ is a non-negative real number.
The set of all functions $f(z)$ which are holomorphic in $T(C^\prime;r)$ and
satisfy the estimate (\ref{eq31}) will be denoted by $\boldsymbol{{\mathscr H}^o_c}$.

\begin{remark}
The space of functions $\boldsymbol{{\mathscr H}^o_c}$ constitutes a generalization
of the space ${\mathfrak A}_{_\omega}^i$ of Sebasti\~ao e Silva~\cite{Tiao1} and the
space $\mona_{_\omega}$ of Hasumi~\cite{Hasumi} to arbitrary tubular radial domains
in $\oC^n$.
\end{remark}

\begin{lemma}[\cite{Carmi1}, \cite{DanHenri}]
Let $C$ be an open convex cone, and let $C^\prime \Subset C$. Let
$h(\xi)=e^{k|\xi|}g(\xi)$, $\xi \in \oR^n$, be a function with support
in $C^*$, where $g(\xi)$ is a bounded continuous function on $\oR^n$.
Let $y$ be an arbitrary but fixed point of $\bigl(C^\prime \setminus
\bigl(C^\prime \cap B[0;r]\bigr)\bigr)$. Then $e^{-\langle \xi,y \rangle}
h(\xi) \in L^2$, as a function of $\xi \in \oR^n$.
\label{lemma0}
\end{lemma}

\begin{definition}
We denote by $H^\prime_{C^*}(\oR^n;O)$ the subspace of $H^\prime(\oR^n;O)$
of distributions of exponential growth with support in the cone $C^*$:
\begin{equation}
H^\prime_{C^*}(\oR^n;O)=\Bigl\{V \in H^\prime(\oR^n;O) \mid
\supp(V) \subseteq C^* \Bigr\}\,\,. 
\label{eq31'} 
\end{equation}
\label{Def1}
\end{definition}

\begin{lemma}[\cite{Carmi1}, \cite{DanHenri}]
Let $C$ be an open convex cone, and let $C^\prime \Subset C$.
Let $V=D^\gamma_\xi[e^{h_K(\xi)}g(\xi)]$, where $g(\xi)$ is a bounded continuous
function on $\oR^n$ and $h_K(\xi)=k|\xi|$ for a convex compact set
$K=\bigl[-k,k\bigr]^n$. Let $V \in H^\prime_{C^*}(\oR^n;O)$. Then $f(z)=(2\pi)^{-n}
\bigl\langle V,e^{-i\langle \xi,z \rangle}\bigr\rangle$ is an element of
$\boldsymbol{{\mathscr H}^o_c}$.
\label{lemma1}
\end{lemma}

%%%%%%%%%%%%%%%%%%%%%%%%%%%%%%%%%%%%%%%%%%%%%%%%%%%%%%%%%%%%%%%%%%%%%%%%%%%%%%%%%
\section{The Space of Holomorphic Functions $\boldsymbol{{\mathscr H}^{*\,o}_c}$}
\label{SecTheo2}
%\hspace*{\parindent}
%%%%%%%%%%%%%%%%%%%%%%%%%%%%%%%%%%%%%%%%%%%%%%%%%%%%%%%%%%%%%%%%%%%%%%%%%%%%%%%%%
We now shall introduce another space of holomorphic functions whose elements are
analytic in a domain $T(C^\prime)$ which is larger than $T(C^\prime;r)$ and has
boundary values in $\oR^n$. The boundary values so obtained are of importance in the
representation of vacuum expectation values in the case of a quantum field theory
in non-commutatives spacetimes~\cite{DZH}.

Let $C$ be a proper open convex cone, and let $C^\prime \Subset C$.
Let $B(0;r)$ denote an {\bf open} ball of the origin in $\oR^n$ of radius
$r$, where $r$ is an arbitrary positive real number. Denote $T(C^\prime;r)=
\oR^n+i\bigl(C^\prime \setminus \bigl(C^\prime \cap B(0;r)\bigr)\bigr)$. Throughout
this section, we consider functions $f(z)$ which are holomorphic in
$T(C^\prime)=\oR^n+iC^\prime$ and which satisfy the estimate (\ref{eq31}),
with $B[0;r]$ replaced by $B(0;r)$. We denote this space by
$\boldsymbol{{\mathscr H}^{*\,o}_c}$. We note that $\boldsymbol{{\mathscr H}^{*\,o}_c}
\subset \boldsymbol{{\mathscr H}^{o}_c}$ for any open convex cone $C$. Put
${\mathscr U}_c=\boldsymbol{{\mathscr H}^{*\,o}_c}/\boldsymbol{\Pi}$, that is,
${\mathscr U}_c$ is the quotient space of $\boldsymbol{{\mathscr H}^{*\,o}_c}$
by set of pseudo-polynomials $\boldsymbol{\Pi}$.

\begin{definition}
The set ${\mathscr U}_c$ is the subspace of the tempered ultrahyperfunctions generated
by $\boldsymbol{{\mathscr H}^{*\,o}_c}$ corresponding to a proper open convex cone
$C \subset \oR^n$.
\end{definition}

The following theorems will be important to us in the proof of edge of the wedge
theorem.

\begin{theorem}
Let $C$ be an open convex cone, and let $C^\prime \Subset C$.
Let $V=D^\gamma_\xi h(\xi)$, where $h(\xi)=e^{h_K(\xi)}g(\xi)$ with $g(\xi)$ being a
bounded continuous function on $\oR^n$ and $h_K(\xi)=k|\xi|$ for a convex compact set
$K=\bigl[-k,k\bigr]^n$. Let $V \in H^\prime_{C^*}(\oR^n;O)$. Then

\,\,\,$(i)\quad f(z)=(2\pi)^{-n}\bigl\langle V,e^{-i\langle \xi,z \rangle}\bigr\rangle$
is an element of $\boldsymbol{{\mathscr H}^{*\,o}_c}$,

\,\,\,$(ii)\quad \bigl\{f(z) \mid y={\rm Im}\,z \in C^\prime \Subset C, |y| \leq Q\bigr\}$
is a strongly bounded set in ${\mathfrak H}^\prime(T(O))$, where $Q$ is an arbitrarily but
fixed positive real number,

\,\,\,$(iii)\quad f(z)\!\rightarrow\!{\mathscr F}^{-1}[V] \in {\mathfrak H}^\prime(T(O)\!)$
in the strong $($and weak$)$ topology of ${\mathfrak H}^\prime(T(O))$ as $y={\rm Im}\,z
\rightarrow 0$, $y \in C^\prime \Subset C$.
\label{theorem1}
\end{theorem}

\begin{theorem}
Let $f(z) \in \boldsymbol{{\mathscr H}^{*\,o}_c}$, where $C$ is an open convex cone.
Then the distribution $V \in H^\prime_{C^*}(\oR^n;O)$ has a uniquely
determined inverse Fourier-Laplace transform $f(z)=(2\pi)^{-n}
\bigl\langle V,e^{-i\langle \xi,z \rangle} \bigr\rangle$ which is holomorphic in
$T(C^\prime)$ and satisfies the estimate (\ref{eq31}), with $B[0;r]$ replaced by
$B(0;r)$.
\label{PWSTheo} 
\end{theorem}

The Theorem \ref{theorem1} shows that functions in
$\boldsymbol{{\mathscr H}^{*\,o}_c}$ have distributional boundary values in
${\mathfrak H}^\prime(T(O))$. Further, it shows that functions in
$\boldsymbol{{\mathscr H}^{*\,o}_c}$ satisfy a strong boundedness property in
${\mathfrak H}^\prime(T(O))$. On the other hand, the Theorem \ref{PWSTheo} shows that
the functions $f(z) \in \boldsymbol{{\mathscr H}^{*\,o}_c}$ can be recovered as
the (inverse) Fourier-Laplace transform of the constructed distribution $V \in H^\prime_{C^*}
(\oR^n;O)$. This result is a version of the Paley-Wiener-Schwartz theorem
in the tempered ultrahyperfunction set-up.

\begin{remark}
It is important to note that in Theorems \ref{theorem1} and \ref{PWSTheo}
we are considering the inverse Fourier-Laplace transform $f(z)=(2\pi)^{-n}
\bigl\langle V,e^{-i\langle \xi,z \rangle} \bigr\rangle$, in contrast to the
Fourier-Laplace transform used in the Refs.~\cite{Carmi1}, \cite{Carmi2}.
\label{remark2}
\end{remark}

\begin{proof}[Sketch of Proof of Theorem \ref{theorem1}]
In order to prove $(i)$, we can proceed as in the proof of~\cite[Lemma 2]{DanHenri}
and obtain the equality
\begin{equation}
f(z)=(2\pi)^{-n}\bigl\langle V,e^{-i\langle \xi,z \rangle}\bigr\rangle
\,\,,\quad z \in T(C^\prime;r)\,\,,
\label{eq311}
\end{equation}
with $B[0;r]$ replaced by $B(0;r)$ in the estimate (\ref{eq31}). The equality
(\ref{eq311}) holds pointwise for arbitrary compact subcones $C^\prime$ of $C$
and for arbitrary $r > 0$. Since $C$ is open, then for any $y \in C$ there is
a compact subcone $C^\prime$ of $C$ and a $r > 0$ such that $y \in \bigl(C^\prime
\setminus \bigl(C^\prime \cap B(0;r)\bigr)\bigr)$. Hence any $z \in T(C^\prime)$ is
in $T(C^\prime;r)$ for some $C^\prime \subset C$ and some $r > 0$. Thus we can conclude
that $(i)$ is obtained from (\ref{eq311}). The proofs of $(ii)$ and $(iii)$ are
similar to the proofs of the Equations (35) and (36) in Ref.~\cite[Theorem 3]{Carmi1}. 
\end{proof}

\begin{proof}[Proof of Theorem \ref{PWSTheo}]
Consider
\begin{equation}
h_y(\xi)=\int_{\oR^n}\frac{f(z)}{P(iz)}\,\,
e^{i\langle \xi,z \rangle}d^nx\,\,,\quad z \in T(C^\prime;r)\,\,, 
\label{eq41p}
\end{equation}
with $h_y(\xi)=e^{k|\xi|}g_y(\xi)$, where $g(\xi)$ is a bounded continuous
function on $\oR^n$, and $P(iz)=(-i)^{|\gamma|} z^\gamma$. By
hypothesis $f(z) \in \boldsymbol{{\mathscr H}^{*\,o}_c}$ and satisfies
(\ref{eq31}), with $B[0;r]$ replaced by $B(0;r)$. For this reason, for an
$n$-tuple $\gamma=(\gamma_1,\ldots,\gamma_n)$ of non-negative integers
conveniently chosen, we obtain
\begin{equation}
\Bigl|\frac{f(z)}{P(iz)}\Bigr|\leq
{\boldsymbol{\sf M}}(C^\prime)(1+|z|)^{-n-\varepsilon} e^{h_{c^*}(y)}\,\,,
\label{eq42} 
\end{equation}
where $n$ is the dimension and $\varepsilon$ is any fixed positive real
number. This implies that the function $h_y(\xi)$ exists and is a continuous function
of $\xi$. Further, by using arguments paralleling the analysis in~\cite[p.225]{Vlad}
and the Cauchy-Poincar\'e Theorem~\cite[p.198]{Vlad}, we can show that the
function $h_y(\xi)$ is independent of $y={\rm Im}\,z$. Therefore, we denote
the function $h_y(\xi)$ by $h(\xi)$.

From (\ref{eq42}) we have that $f(z)/P(iz) \in L^2$ as a function of
$x={\rm Re}\,z \in \oR^n$, $y \in C^\prime \setminus \bigl(C^\prime
\cap B(0;r)\bigr)$. Hence, from (\ref{eq41p}) and the Plancherel theorem
we have that $e^{-\langle \xi,y \rangle}h(\xi) \in L^2$ as a function of
$\xi \in \oR^n$, and  
\begin{equation}
\frac{f(z)}{P(iz)}={\mathscr F}^{-1}\bigl[e^{-\langle \xi,y \rangle}
h(\xi)\bigr](x)\,\,,\quad z \in T(C^\prime;r)\,\,,
\label{eq43} 
\end{equation}
where the inverse Fourier transform is in the $L^2$ sense. Here,
Parseval's equation holds:
\begin{equation}
(2\pi)^{-n}\int_{\oR^n}\Bigl|e^{-\langle \xi,y \rangle}h(\xi)\Bigr|^2d^n\xi=
\int_{\oR^n}\Bigl|\frac{f(z)}{P(iz)}\Bigr|^2
d^nx\,\,. 
\label{eq44}
\end{equation}
In this case for the Eq.(\ref{eq43}) to be true $\xi$ must belong to the
open half-space $\bigl\{\xi \in C^* \mid \langle \xi,y \rangle < 0\bigr\}$, for
$y \in C^\prime \setminus \bigl(C^\prime \cap B(0;r)\bigr)$, since by hypothesis
$f(z) \in \boldsymbol{{\mathscr H}^{*\,o}_c}$. Then there is $\delta(C^\prime)$
such that for $y \in C^\prime \setminus \bigl(C^\prime \cap B(0;r)\bigr)$ implies
$\langle \xi,y \rangle \leq -\delta(C^\prime) |\xi||y|$. This justifies the
negative sign in (\ref{eq43}) (see Remark \ref{remark2}).

Now, if $h(\xi) \in H^\prime_{C^*}(\oR^n;O)$, then $V=D^\gamma_\xi h(\xi)
\in H^\prime_{C^*}(\oR^n;O)$. Since $C^*$ is a regular set~\cite[pp.98, 99]{Sch},
thus $\supp(h)=\supp(V)$. By Theorem \ref{theorem1} $\bigl\langle V,e^{-i\langle
\xi,z \rangle} \bigr\rangle$ exists as a holomorphic function of $z \in T(C^\prime)$
and satisfies the estimate (\ref{eq31}), with $B[0;r]$ replaced by $B(0;r)$.
A simple calculation yields
\begin{equation}
(2\pi)^{-n}\bigl\langle V,e^{-i\langle \xi,z \rangle} \bigr\rangle=
P(iz){\mathscr F}^{-1}\bigl[e^{-\langle \xi,y \rangle}
h(\xi)\bigr](x)\,\,\quad z \in T(C^\prime;r)\,\,.
\label{eq45} 
\end{equation}
In view of Lemma \ref{lemma0}, the inverse Fourier transform can be interpreted
in $L^2$ sense. Combining (\ref{eq43}) and (\ref{eq45}), we have
$f(z)=(2\pi)^{-n}\bigl\langle V,e^{-i\langle \xi,z \rangle} \bigr\rangle$,
for $z \in T(C^\prime;r)$. Since $r > 0$ is arbitrary, this equality holds for
each $z \in T(C^\prime)$. The uniqueness follows from the isomorphism of the dual
Fourier transform, according to Proposition \ref{Propo1}. This completes the proof
of the theorem.  
\end{proof}

%%%%%%%%%%%%%%%%%%%%%%%%%%%%%%%%%%%%%%%%%%%%%%%%%%%%%%%%%%%%%%%%%%%%%%%%%%%%%%%%%
\section{Edge of the Wedge Theorem}
\label{SecTheo3}
%\hspace*{\parindent}
%%%%%%%%%%%%%%%%%%%%%%%%%%%%%%%%%%%%%%%%%%%%%%%%%%%%%%%%%%%%%%%%%%%%%%%%%%%%%%%%%
In what follows, we formulate a version of the edge of the wedge theorem for the
space of the tempered ultrahyperfunctions in its simplest form: the common analytic
continuation of two functions $f_1(z)$ and $f_2(z)$
holomorphic respectively in the two tubes $\oR^n+iC_j$, $j=1,2$, where
each $C_j$ is an open convex cone.

\begin{theorem}[Edge of the Wedge Theorem]
Let $C$ be an open cone of the form $C=C_1 \cup C_2$, where each $C_j$, $j=1,2$,
is a proper open convex cone. Denote by $\boldsymbol{ch}(C)$ the convex hull of the
cone $C$. Assume that the distributional boundary values of two holomorphic functions
$f_j(z) \in \boldsymbol{{\mathscr H}^{*\,o}_{c_j}}$ $(j=1,2)$ agree, that
is, $U=BV(f_1(z))=BV(f_2(z))$, where $U \in {\mathfrak H}^\prime(T(O))$ in accordance
with the Theorem \ref{theorem1}. Then there exists $F(z) \in
\boldsymbol{{\mathscr H}^{o}_{{\boldsymbol{ch}(C)}}}$
such that $F(z)=f_j(z)$ on the domain of definition of each $f_j(z)$, $j=1,2$.
\label{EWTheo}  
\end{theorem}

\begin{proof}
By hypothesis $BV(f_1(z))=BV(f_2(z))$ in ${\mathfrak H}^\prime(T(O))$, and we call
this common value $U$. By Theorem \ref{theorem1}, we have that $BV(f_j(z))={\mathscr F}^{-1}[V_j]$,
$j=1,2$. On the other hand, this implies that $V_j={\mathscr F}[BV(f_j(z))]$. But, according to
Theorem \ref{PWSTheo} there exists a unique $V_j \in H^\prime_{C_j^*}(\oR^n;O)$, $j=1,2$,
such that $f_j(z)=(2\pi)^{-n}\bigl\langle V_j,e^{-i\langle \xi,z \rangle} \bigr\rangle$.
Using these facts we have that $V_1=V_2$ in $H^\prime_{C^*}(\oR^n;O)$. We call this common
value $V$ and thus have $U={\mathscr F}^{-1}[V]$. By Theorem 2 in~\cite{Carmi3},
$\supp(V) \subseteq \bigl\{\xi \in \oR^n \mid \langle \xi,y \rangle \geq 0,\forall\,\,y
\in \boldsymbol{ch}(C)\bigr\}$, then by Definition \ref{Def1} $V \in
H^\prime_{(\boldsymbol{ch}(C))^*}(\oR^n;O)$.
 
We now put
\begin{equation}
F(z)=(2\pi)^{-n}\bigl\langle V,e^{-i\langle \xi,z \rangle} \bigr\rangle\,\,,
\quad z \in T(\boldsymbol{ch}(C))=\oR^n+i\,\boldsymbol{ch}(C)\,\,.
\label{eq01}
\end{equation}
with $V \in H^\prime_{(\boldsymbol{ch}(C))^*}(\oR^n;O)$.
Since $\boldsymbol{ch}(C)$ is an open convex cone, we have by exactly the proof
of~\cite[Lemma 2]{DanHenri} that $F(z) \in \boldsymbol{{\mathscr H}^{o}_{{\boldsymbol{ch}(C)}}}$.
Further, using the fact that $V_1=V_2=V$, from Theorem \ref{theorem1} we have that
\begin{equation}
f_j(z)=(2\pi)^{-n}\bigl\langle V_j,e^{-i\langle \xi,z \rangle} \bigr\rangle=
(2\pi)^{-n}\bigl\langle V,e^{-i\langle \xi,z \rangle} \bigr\rangle\,\,,
\quad z \in T(C^\prime_j)\,\,.
\label{eq02}
\end{equation}
Thus combining (\ref{eq01}) and (\ref{eq02}) we have that $F(z)$ coincides with $f_j(z)$,
$j=1,2$, on the domain of definition of each $f_j(z)$. 
\end{proof}

\begin{corollary}
Suppose that the hypotheses of Theorem \ref{EWTheo} hold with $C_1$ and $C_2$ opposite to
each other. Then $F(z)$ is a polynomial in $z \in \oC^n$.
\label{Corollary1} 
\end{corollary}

\begin{proof}
Similar to the proof of~\cite[Corollary 1]{Carmi3}.
\end{proof}

The following theorem is an immediate consequence of the edge of the wedge theorem and
reflects a of the most important principle governing the behaviour of analytic functions,
that is, the determination of a function by its values on a non-empty open real set.

\begin{theorem}
Let $C$ be some open convex cone. Let $f(z) \in \boldsymbol{{\mathscr H}^{*\,o}_{c}}$.
If the distributional boundary value $BV(f(z))$ of $f(z)$ in the sense of tempered
ultrahyperfunctions vanishes, then the function $f(z)$ itself vanishes.
\label{SpCase}
\end{theorem}

\begin{proof}
Define $g(x+iy)=\overline{f(x-iy)}$. The function $g(z)$ is holomorphic in
$\overline{T(C^\prime)}=\oR^n-iC^\prime$, satisfies (\ref{eq31}), with $B[0;r]$ replaced
by $B(0;r)$, in $\overline{T(C^\prime;r)}=\oR^n-i\bigl(C^\prime \setminus \bigl(C^\prime
\cap B(0;r) \bigr)\bigr)$ and approaches $0$ as $y \rightarrow 0$. Thus we can apply the
edge of the wedge theorem to $f$ and $g$. Since $\boldsymbol{ch}(C \cup (-C))=\oR^n$,
then by Corollary \ref{Corollary1} $F(z)$, the common analytic continuation of $f$ and
$g$, is a polynomial in $z \in \oC^n$. But, by hypothesis $BV(F(z))$ vanishes as a
distribution and therefore as a function together with $f(z)$ identically.
\end{proof}

%%%%%%%%%%%%%%%%%%%%%%%%%%%%%%%%%%%%%%%%%%%%%%%%%%%%%%%%%%%%%%%%%%%%%%%%%%%%%%%%%
\section{The Martineau's Edge of the Wedge Theorem\\
for Tempered Ultrahyperfunctions}
\label{SecTheo4}
%\hspace*{\parindent}
%%%%%%%%%%%%%%%%%%%%%%%%%%%%%%%%%%%%%%%%%%%%%%%%%%%%%%%%%%%%%%%%%%%%%%%%%%%%%%%%%
The great advance in the theory of the edge of the wedge theorem came with the realization
due to Martineau~\cite{Marti1}, \cite{Marti2}, \cite{Marti3}, who was able to prove its
version for the case involving more than two functions holomorphic respectively in the
tubes $\oR^n+iC_j$, $j=1,\ldots,m$. In what follows, we formulate a version of the Martineau's
edge of the wedge theorem for the space of the tempered ultrahyperfunctions.

\begin{theorem}[Generalized Edge of the Wedge Theorem]
Let $C_1,\ldots,C_m$ be proper open convex cones in $\oR^n$. Given any set of $m$
open convex cones $C^\prime_j$ such that $C^\prime_j \Subset C_j$, $j=1,\ldots,m$,
then the following two properties of a set of $m$ functions $f_j(z) \in
\boldsymbol{{\mathscr H}^{*\,o}_{c_j}}$ $(j=1,\ldots,m)$ are equivalent:
%%%%%%%%%%%%%%%%%%%%%%%%%%%%%%%%%%%%%%%%%%%%%%%%%%%%%%%%%%%%%%%%%%%%%%%%%%%%%%%%%%%%%%%%%
\newcounter{numero}
\setcounter{numero}{0}
\def\Prop{\addtocounter{numero}{1}\item[{$\boldsymbol{\sf P_{\thenumero}-}$}]}
\begin{enumerate}

\Prop The distributional boundary value $U=\sum_{j=1}^{m} BV(f_j(z)) \in
{\mathfrak H}^\prime(T(O))$ vani\-shes identically.

\Prop Denote by $\boldsymbol{ch}(C_j \cup C_k)$ the convex hull of $C_j \cup C_k$.
For each pair of indices $(j,k)$, $1 \leq j,k \leq m$, there is a holomorphic
function $g_{jk}(z) \in \boldsymbol{{\mathscr H}^{o}_{{\boldsymbol{ch}(C_j \cup C_k)}}}$ 
such that $g_{jk}(z)+g_{kj}(z)=0$ for all $j,k=1,\ldots,m$ -- thus $g_{jj}(z)=0,\,\,
\forall\,\,j,k=1,\ldots,m$ -- and such moreover that $f_j(z)=\sum_{k=1}^{m} g_{jk}(z)$
on $T(C_j^\prime)=\oR^n+iC_j^\prime$, for each $j=1,\ldots,m$.
\end{enumerate} 
%%%%%%%%%%%%%%%%%%%%%%%%%%%%%%%%%%%%%%%%%%%%%%%%%%%%%%%%%%%%%%%%%%%%%%%%%%%%%%%%%%%%%%%%
\label{GEWTheo}  
\end{theorem}

For our proof of the Theorem \ref{GEWTheo} we prepare a lemma on the analytic
decomposability of ${\mathfrak H}^\prime(T(O))$. Let $C$ be an open cone of the form
$C=\bigcup_{j=1}^m C_j$, $m < \infty$, where each $C_j$ is an proper open convex cone.
If we write $C^\prime \Subset C$, we mean $C^\prime=\bigcup_{j=1}^m C_j^\prime$ with
$C_j^\prime \Subset C_j$. Furthermore, we define by $C_j^*=\bigl\{\xi \in \oR^n \mid
\langle \xi,x \rangle \geq 0, \forall x \in C_j \bigr\}$ the dual cones of $C_j$, such
that the dual cones $C_j^*$, $j=1,\ldots,m$, have the properties
\begin{equation}
\oR^n \setminus \bigcup_{j=1}^m C_j^*\,\,,
\label{p1}
\end{equation}
and
\begin{equation}
C_j^* \bigcap C_k^*\,\,,j \not= k\,\,,j,k=1,\ldots,m\,\,,
\label{p2}
\end{equation}
are sets of Lebesgue measure zero. Assume that $V \in H^\prime_{C^*}(\oR^n;O)$
can be written as $V=\sum_{j=1}^{m} V_j$, where we define
\begin{equation}
V_j=D^\gamma_\xi[e^{h_K(\xi)}\lambda_j(\xi)g(\xi)]\,\,,
\label{Def2}
\end{equation}
with $\lambda_j(\xi)$ denoting the characteristic function of $C_j^*$, $j=1,\ldots,m$,
$g(\xi)$ being a bounded continuous function on $\oR^n$ and $h_K(\xi)=k|\xi|$ for a
convex compact set $K=\bigl[-k,k\bigr]^n$.

\begin{lemma}
Let $C$ be an open cone of the form $C=\bigcup_{j=1}^m C_j$, $m < \infty$,
where the $C_j$ are proper open convex cones such that $($\ref{p1}$)$ and $($\ref{p2}$)$
are satisfied. Let $U \in {\mathfrak H}^\prime(T(O))$. Then $U=\sum_{j=1}^{m} BV(f_j(z))$,
where each $BV(f_j(z))$ is the strong boun\-da\-ry value in ${\mathfrak H}^\prime(T(O))$
of a function $f_j(z) \in \boldsymbol{{\mathscr H}^{*\,o}_{c_j}}$ and such that each
$BV(f_j(z))\!=\!{\mathscr F}^{-1}[V_j]$, where $V_j \in H^\prime_{C_j^*}(\oR^n;O)$,
$j=1,\ldots,m$.
\end{lemma}

\begin{proof}
This result follows using the same method adopted in the proof of~\cite[Theorem 4]{Carmi1},
by replacing the reference to Theorems 2 and 3 by a reference to the Theorem \ref{theorem1}
of this paper.
\end{proof}
 
\begin{proof}[Proof of Theorem \ref{GEWTheo}]
We begin by proving that $\boldsymbol{\sf P_2} \Rightarrow \boldsymbol{\sf P_1}$.
Assume $g_{jk}(z) \in \boldsymbol{{\mathscr H}^{o}_{{\boldsymbol{ch}(C_j \cup C_k)}}}$.
By hypothesis, we have
\[
f_j(z)=\sum_{k=1}^{m} g_{jk}(z)\,\,,\quad z \in T(C_j^\prime)\,\,. 
\]
Then,
\[
BV(f_j(z))=BV\Bigl(\sum_{k=1}^{m} g_{jk}(z)\Bigr)\,\,,
\quad {\mbox{as $C_j^\prime \ni y \rightarrow 0$}}\,\,. 
\]
Hence,
\begin{align*}
\sum_{j=1}^{m} BV(f_j(z))&=\sum_{j=1}^{m}\Bigl(BV \bigl(\sum_{k=1}^{m}
g_{jk}(z)\bigr)\Bigr)\\[3mm]
&=BV \Bigl(\sum_{j=1}^{m} \sum_{k=1}^{m}
g_{jk}(z)\Bigr) \equiv 0\,\,,
\end{align*}
taking into account the anti-symmetry of the functions $g_{jk}(z)$.

Proof that $\boldsymbol{\sf P_1} \Rightarrow \boldsymbol{\sf P_2}$. If $m=1$,
$\boldsymbol{\sf P_1} \Rightarrow f_1 \equiv 0$ by Theorem \ref{SpCase}. Henceforth
we assume $m \geq 2$. Let $U_j=BV(f_j(z)) \in {\mathfrak H}^\prime(T(O))$ as $C_j^\prime
\ni y \rightarrow 0$. Then there exists $U_{jk} \in {\mathfrak H}^\prime(T(O))$ such
that $U_j=\sum_{k=1}^{m} U_{jk}$ with the restriction that $U_{jk}+U_{kj}=0$. Thus $\sum_{j=1}^{m}
U_j \equiv 0$. By Theorem \ref{theorem1}, we have that $U_j={\mathscr F}^{-1}[V_j]$,
$j=1,\ldots,m$. On the other hand, this implies that $V_j={\mathscr F}[U_j]={\mathscr F}
[\sum_{k=1}^{m} U_{jk}]=\sum_{k=1}^{m}{\mathscr F}[U_{jk}]=\sum_{k=1}^{m} V_{jk}$.
Since $\sum_{j=1}^{m} U_j \equiv 0$, it follows that ${\mathscr F}[\sum_{j=1}^{m} U_j]=
\sum_{j=1}^{m} {\mathscr F}[U_j]=\sum_{j=1}^{m} V_j=\sum_{j=1}^{m}\sum_{k=1}^{m} V_{jk}
\equiv 0$. This implies that $V_{jk}+V_{kj}=0$. According to Theorem \ref{PWSTheo} there
exists a unique $V_j \in H^\prime_{C_j^*}(\oR^n;O)$, $j=1,\ldots,m$, such that
\begin{align*}
f_j(z)=(2\pi)^{-n}\bigl\langle V_j,e^{-i\langle \xi,z \rangle} \bigr\rangle
&=(2\pi)^{-n}\bigl\langle\sum_{k=1}^{m} V_{jk},e^{-i\langle \xi,z \rangle}
\bigr\rangle \\[3mm]
&=\sum_{k=1}^{m} \Bigl((2\pi)^{-n}\bigl\langle V_{jk},e^{-i\langle \xi,z \rangle}
\bigr\rangle \Bigr)\,\,.
\end{align*}
We now put
\begin{equation}
g_{jk}(z)=(2\pi)^{-n}\bigl\langle V_{jk},e^{-i\langle \xi,z \rangle} \bigr\rangle\,\,,
\quad z \in T(\boldsymbol{ch}(C_j \cup C_k))\,\,.
\label{eq011}
\end{equation}
with $V_{jk} \in H^\prime_{(\boldsymbol{ch}(C_j \cup C_k))^*}(\oR^n;O)$ and
$\supp(V_{jk}) \subseteq \bigl\{\xi \in \oR^n \mid \langle \xi,y \rangle \geq 0,\forall\,\,y
\in \boldsymbol{ch}(C_j \cup C_k)\bigr\}$. Since $\boldsymbol{ch}(C_j \cup C_k)$ is an open
convex cone, we again have by exactly the proof of~\cite[Lemma 2]{DanHenri}
that $g_{jk}(z) \in \boldsymbol{{\mathscr H}^{o}_{{\boldsymbol{ch}(C_j \cup C_k)}}}$.
Further, from Theorem \ref{theorem1}, we have that
\begin{align}
f_j(z)&=(2\pi)^{-n}\bigl\langle V_j,e^{-i\langle \xi,z \rangle} \bigr\rangle \nonumber \\[3mm]
&=\sum_{k=1}^{m} \Bigl((2\pi)^{-n}\bigl\langle V_{jk},e^{-i\langle \xi,z \rangle}
\bigr\rangle \Bigr)\,\,,\quad z \in T(C^\prime_j)\,\,.
\label{eq022}
\end{align}
Thus combining (\ref{eq011}) and (\ref{eq022}) we have that $f_j(z)=\sum_{k=1}^{m} g_{jk}(z)$
on $T(C_j^\prime)$ with the restriction that $g_{jk}+g_{kj}=0$. This completes the proof.
\end{proof}

%%%%%%%%%%%%%%%%%%%%%%%%%%%%%%%%%%%%%%%%%%%%%%%%%%%%%%%%%%%%%%%%%%%%%%%%%%%
\section*{Acknowledgments}
%\hspace*{\parindent}
%%%%%%%%%%%%%%%%%%%%%%%%%%%%%%%%%%%%%%%%%%%%%%%%%%%%%%%%%%%%%%%%%%%%%%%%%%%%
The author would like to express his gratitude to Afr\^anio R. Pereira, Winder A.M. Melo
and to the Departament of Physics of the Universidade Federal de Vi\c cosa (UFV)
for the opportunity of serving as Visiting Researcher. The author thanks to the referee
for her/his remarks which improved the presentation of the results.

%%%%%%%%%%%%%%%%%%%%%%%%%%%%%%%%%%%%%%%%%%%%%%%%%%%%%%%%%%%%%%%%%%%%%%%%%%%%%

%%%%%%%%%%%%%%%%%%%%%%%%%%%%%%%%%%%%%%%%%%%%%%%%%%%%%%%%%%%%%%%%%%%%%%%%%%%%%%%%%%%%%

%%%%%%%%%%%%%%%%%%%%%%%%%%%%%%%%%%%%%%%%%%%%%%%%%%%%%%%%%%%%%%%%%%%%%%%%%%%%%%
\end{document}